\documentclass{elsart}
\usepackage{amsmath,amssymb}
\usepackage{array}
\usepackage{epsf}
\numberwithin{equation}{section}

\begin{document}
\newcommand{\derl}{\partial_L}
\newcommand{\derr}{\partial_R}

\def\tr{{\rm tr}\,}
\def\be{\begin{equation}}
\def\ee{\end{equation}}
\def\bea{\begin{eqnarray}}
\def\eea{\end{eqnarray}}
\def\wc{\bar{w}}
\def\ld{\ldots}
\def\ie{i.e.}
\def\C{{\mathbb C}}
\def\Z{{\mathbb Z}}
\def\R{{\mathbb R}}
\def\cf{{\mathcal{F}}}
\def\cs{{\mathcal{S}}}
\def\cl{{\mathcal{L}}}
\def\cx{{\mathcal{X}}}
\def\fd{{f^\dagger}}
\def\fdo{{f^{0\dagger}}}
\def\1{{\bf 1}}

\begin{frontmatter}
\title{Description of surfaces associated with $\C P^{N-1}$ sigma models on Minkowski space}

\author[AMG]{
A.~M. Grundland\thanksref{AMGmail}}
and 
\author[LS]{L. \v{S}nobl\thanksref{LSmail}}
\address[AMG]{
Centre de Recherches Math{\'e}matiques, Universit{\'e} de Montr{\'e}al,
C. P. 6128, Succ.\ Centre-ville, Montr{\'e}al, (QC) H3C 3J7 and
Universit\'{e} du Qu\'{e}bec, Trois-Rivi\`{e}res CP500 (QC) G9A 5H7, Canada }
\address[LS]{
Centre de Recherches Math{\'e}matiques, Universit{\'e} de Montr{\'e}al, 
C. P. 6128, Succ.\ Centre-ville, Montr{\'e}al, (QC) H3C 3J7, Canada and
Faculty of Nuclear Sciences and Physical Engineering, 
Czech Technical University, B\v rehov\'a 7, 115 19 Prague~1, Czech Republic} 
\thanks[AMGmail]{Email address: grundlan@crm.umontreal.ca}
\thanks[LSmail]{Email address: Libor.Snobl@fjfi.cvut.cz}

\begin{abstract}
The objective of this paper is to construct and investigate smooth orientable surfaces in $\R^{N^2-1}$ by analytical methods.
The structural equations of surfaces 
in connection with $\C P^{N-1}$ sigma models on Minkowski space 
are studied in detail. This is carried out using moving frames adapted to surfaces immersed in the $su(N)$ algebra. The first and 
second fundamental forms of this surface as well as the relations between them as expressed in the 
Gauss--Weingarten and Gauss--Codazzi--Ricci equations are found. 
The Gaussian curvature, the mean curvature vector and the Willmore functional expressed in terms of 
a solution of $\C P^{N-1}$ sigma model are obtained. 
An example of a surface associated with the $\C P^1$ model is included as 
an illustration of the theoretical results. 
\end{abstract}
\end{frontmatter}

\section{Introduction}

Over the last two decades surfaces immersed in multi--dimensional spaces have been extensively researched in 
connection with integrable systems  
(for a review see \cite{Rog} and the references therein). The motivation for this activity came largely 
from applications in various branches of physical, biological and chemical sciences as well as from engineering. 
The progress in  the analytical treatment of surfaces obtained from nonlinear differential equations has been rapid and resulted 
in many new techniques and theoretical approaches. Some of the most interesting developments have been in the study of surfaces 
immersed in Lie algebras, using techniques of completely integrable systems \cite{Bob,FokG,HelF}. These surfaces are characterized 
by fundamental forms whose coefficients satisfy the Gauss--Weingarten and the Gauss--Codazzi--Ricci equations.

In this work we apply a group--theoretical approach to surfaces associated with $\C P^{N-1}$ sigma models. 
This line of investigation was initiated in \cite{GZ} where it was shown that two--dimensional 
constant mean curvature surfaces in three-- and eight--dimensional
spaces are associated with the $\C P^{1}$ and $\C P^{2}$ sigma models defined on Euclidean spaces. Further, it was demonstrated 
in \cite{GSZ}
that any surface described by $\C P^{N-1}$ models on Euclidean space can be constructed by a choice of a moving frame based on $su(N)$ 
algebra representation parametrized by a corresponding solution of the model. 
This has been a significant result since surfaces immersed in Lie algebras are known to show up in many physical systems (see e.g. 
\cite{Nel,Cha}). Our objective in this paper is to extend this approach to the case of $\C P^{N-1}$ sigma models defined on 
Minkowski space. 
To this end we have devised a new technique for construction of a moving frame; the properties of surfaces obtained
in this way, e.g. curvatures, turned out to be significantly different from the ones in the case of sigma models on Euclidean space.

The use of sigma models in mathematical physics has encompassed predominantly models defined on Euclidean spaces, 
since a great number of physical systems can be reduced to these models. 
However, in recent literature we find an increasing number of examples when reductions lead to sigma models on Minkowski 
space and the need for description of surfaces related to these models is certainly there. 
One such example is the string theory in which sigma models on spacetime and their supersymmetric extensions play a crucial role.
Classical configuration of strings can be described by common solutions of the Nambu--Goto--Polyakov action and a system of 
Dirac type equations intimately connected to $\C P^{N-1}$ models \cite{Kon,Kon1}.
Other relevant applications of recent interest are in the areas of statistical physics (for example reduction of self--dual 
Yang--Mills equations to the Ernst model \cite{Abl,David}), phase transitions (e.g. dynamics of vortex sheets, 
growth of crystals, surface waves etc.\cite{Nel,Cha}) and the theory of fluid membranes \cite{Ou,Saf}. In this last example it is known 
that the free energy per molecule is determined by two invariants (the mean and Gauss curvatures) of a surface associated 
with particular class of solutions of $\ C P^{N-1}$ sigma model defined on Minkowski space.

\section{$\C P^{N-1}$ sigma models and their Euler--Lagrange equations}

As a starting point let us present some basic formulae and notation for $\C P^{N-1}$ sigma models defined on Minkowski space.

The points  of the complex coordinate space $\C^{N}$ will be
denoted by $z=(z_1,\ldots,z_{N})$ and the hermitian inner product
in $\C^{N}$ by 
\be\label{1.1} \langle z,w \rangle=z^\dagger
w=\sum_{j=1}^{N}\bar{z}_j w_j\,.
\ee
The complex projective space $\C P^{N-1}$ is defined as a set of 1--dimensional subspaces in $\C^N$. The manifold structure on it is 
defined by an open covering 
$$U_{k} = \{ [z] \; | z \in \C^N, z_k \neq 0 \}, \; k=1,\ldots,N, $$
where $[z]={\rm span} \{ z \}$, 
and coordinate maps 
$$ \varphi_k: U_k \rightarrow \C^{N-1}, \; \; 
\varphi_k(z) =\left( \frac{z_1}{z_k}, \ldots,{\frac{z_{k-1}}{z_k}},{\frac{z_{k+1}}{z_k}}, \ldots,\frac{z_N}{z_k} \right).$$

Let $\xi^1$, $\xi^2$ be the standard Minkowski coordinates in $\R^2$, with the metric
$$ {\rm d}s^2 =({\rm d}\xi^1)^2-({\rm d}\xi^2)^2. $$ 
In what follows we suppose that
$\xi_L=\xi^1+\xi^2 ,\, \xi_R=\xi^1-\xi^2$ are the light--cone coordinates in $\R^2$, i.e.     
\be\label{Minkmetr}
 {\rm d}s^2= {\rm d}\xi_L {\rm d}\xi_R. 
\ee
We shall denote by $\derl$ and $\derr$
the derivatives with respect to $\xi_L$ and $\xi_R$, respectively, \ie
$$\derl=\frac{1}{2}\left(\partial_{\xi^1}+\partial_{\xi^2}\right)\,,\quad
\derr=\frac{1}{2}\left(\partial_{\xi^1}-\partial_{\xi^2}\right)\,.$$ 
Let us assume that $\Omega$ is an open, connected and simply connected subset in
$\R^2$ with Minkowski metric (\ref{Minkmetr}).
In the study of $\C P\sp{N-1}$ sigma models we are interested in maps $[z]:\Omega\rightarrow \C P\sp{N-1}$
 which are stationary points of the action functional (see e.g. \cite{Z}) 
\be\label{cpmodel}
{\cs} = \frac{1}{4} \int_\Omega (D_\mu z)^{\dagger} (D^\mu z) {\rm d}\xi^1 {\rm d}\xi^2, \qquad z^\dagger . z=1.
\ee
The covariant derivatives $D_\mu$ act on $z: \Omega \rightarrow \C^{N}$ according to the formula
\be
D_\mu z = \partial_{\mu} z - (z^{\dagger} . \partial_\mu z) z, \qquad \partial_{\mu} \equiv \partial_{\xi^\mu}, \ \mu=1,2
\ee
and ensure that the action depends only on $[z]:\Omega\rightarrow \C P\sp{N-1}$ and not on the choice of a representative 
of the class $[z]$. 
Thus the map $[z]$ is determined as a solution of the 
Euler--Lagrange equations defined by the action (\ref{cpmodel}). Writing 
\be
z = \frac{f}{|f|}, \; \; |f|= \sqrt{\fd f}
\ee
one can present the action functional (\ref{cpmodel}) also in the form
\be\label{modact}
{\cs} = \int_\Omega {\cl} {\rm d}\xi_L {\rm d}\xi_R = \int_\Omega \frac{1}{4|f|^2} \left( \derl f^\dagger P \derr f + \derr f^\dagger P \derl f \right)  
{\rm d}\xi_L {\rm d}\xi_R,
\ee
where the $N \times N$ matrix
\be\label{defP}
P = \1-\frac{1}{|f|^2} f \otimes f^\dagger
\ee
is an orthogonal projector on $\C^{N}$, i.e. $ P^2=P, \; P^\dagger=P$.

It is useful to recall that the action (\ref{modact}) has the local (gauge) $U(1) \times \R$ symmetry 
\be\label{U1}
f \rightarrow {\rm e}^{i \alpha(\xi_L,\xi_R)+\beta(\xi_L,\xi_R)} f, \; \alpha(\xi_L,\xi_R),\beta(\xi_L,\xi_R):\Omega \rightarrow \R
\ee
corresponding to the fact that the model is defined on $\C P^{N-1}$. In addition, the action (\ref{modact}) has 
the $U(N)$ global symmetry
\be\label{UN}
 f \rightarrow \Phi f, \; \Phi \in U(N).\ee 
It is also invariant under the conformal transformations 
\be\label{Poincare}
 \xi_L \rightarrow \alpha (\xi_L), \; \xi_R \rightarrow \beta ( \xi_R ),  
\ee
where $\alpha,\beta: \R \rightarrow \R$ are arbitrary 1--to--1 maps such that $\derl \alpha(\xi_L)\neq 0, \ \derr \beta(\xi_R) \neq 0$,
as well as under the parity transformation 
\be\label{parity}
 \xi_L \rightarrow \xi_R, \; \xi_R \rightarrow \xi_L .
\ee
Let us note that the invariance properties (\ref{U1})--(\ref{parity}) are naturally reproduced on the level of Euler--Lagrange equations.

Computing the extremals of the action (\ref{modact}), one finds 
the Euler--Lagrange equations in terms of $f$ 
\be\label{eqnmotN}
 P \  \{ \derl \derr f - \frac{1}{(f^\dagger f)} \left( (f^\dagger \derr f) \derl f + (f^\dagger \derl f) \derr f  \right)  \} 
 = 0.
\ee
They can be also expressed in the matrix form
\be\label{eqnmotmat}
[\derl \derr P,P] =0
\ee
or in the form of a conservation law
\be\label{eqnmotmatcons}
\derl [\derr P, P] + \derr [\derl P,P] =0.
\ee
By explicit calculation one can check that the real--valued functions
\be\label{currents}
J_L = \frac{1}{f^\dagger f} \derl f^\dagger P \derl f, \; \; J_R = \frac{1}{f^\dagger f} \derr f^\dagger P \derr f
\ee
satisfy
\be\label{conslawN}
 \derl J_R = \derr J_L =0 
\ee
for any solution $f$ of the Euler--Lagrange equations (\ref{eqnmotN}). The functions $J_L,J_R$ are invariant under local 
$U(1) \times \R$ and global $U(N)$ transformations (\ref{U1}) and (\ref{UN}).

\section{Surfaces obtained from $\C P^{N-1}$ sigma model}\label{Nsurf}

Let us now discuss the analytical description of a two--dimensional smooth orientable surface $\cf$ immersed in the $su(N)$ algebra, 
associated with $\C P^{N-1}$ sigma model (\ref{eqnmotN}).
We shall construct an exact $su(N)$--valued 1--form whose ``potential'' 0--form defines the surface $\cf$. 
Next, we shall investigate the geometric characteristics of the surface $\cf$. 

Let us introduce a scalar product
$$ ( A , B ) = -\frac{1}{2} \tr AB $$
on $su(N)$ and identify the $(N^2-1)$--dimensional Euclidean space with the $su(N)$ algebra
$$ \R^{N^2-1} \simeq su(N). $$
We denote 
\be M_L = [\derl P,P], \qquad M_R = [\derr P,P], \ee
or, equivalently, using (\ref{defP})
\be\label{explM}
M_{D} = \frac{1}{\fd f} ( P \partial_{D} f \otimes \fd - f \otimes \partial_{D} \fd P) \; \in su(N), \; D=L,R.
\ee
It follows from (\ref{eqnmotmatcons}) that if $f$ is a solution of the Euler--Lagrange equations (\ref{eqnmotN}) then
\be\label{eqnmotmatcons2}
\derl M_R+ \derr M_L =0 .
\ee
Therefore we can identify tangent vectors to the surface $\cf$ with the matrices $M_L$ and $M_R$, as follows
\be\label{dlxdrx}
 {\cx}_L=M_L, \; \; {\cx}_R=-M_R 
\ee
Equation (\ref{eqnmotmatcons2}) implies there exists a closed $su(N)$--valued 1--form on $\Omega$
$$ {\cx} = {\cx}_L {\rm d} \xi_L + {\cx}_R {\rm d} \xi_R, \; \;  {\rm d} {\cx} =0 .$$
Because ${\cx}$ is closed and $\Omega$ is connected and simply connected, ${\cx}$ is also exact. In other words,
there exists a  well--defined 
$su(N)$--valued function $X$ on $\Omega$ such that ${\cx} = {\rm d} X$. 
The matrix function $X$ is unique up to addition of any constant element of $su(N)$ and 
we identify the elements of $X$ with the coordinates of the sought--after surface $\cf$ in $\R^{N^2-1}$. 
Consequently, we get 
\be\label{cldx}
\derl X= {\cx}_L,\ \derr X= {\cx}_R.
\ee
The map $X$ is called the Weierstrass formula for immersion. 
In practice, the surface $\cf$ is found by integration
\be\label{surfaceN}
\cf: \;  X(\xi_L,\xi_R) = \int_{\gamma(\xi_L,\xi_R)} {\cx}
\ee
along any curve $\gamma(\xi_L,\xi_R)$ in $\Omega$ connecting the point $(\xi_L,\xi_R)\in \Omega$ with an arbitrary chosen
point $(\xi_L^0,\xi_R^0)\in \Omega$.

By computation of traces of ${\cx}_{B}.{\cx}_{D}, \ B,D = L,R$ we immediately find the components of the induced metric on the surface 
${\cf}$
\be\label{metric} 
G  =  \left( \begin{array}{cc} G_{LL}, & G_{LR} \\ G_{LR}, & G_{RR} \end{array}  \right) = 
\left( \begin{array}{cc} J_L  & - \Re \left( \frac{\derr \fd P \derl f}{\fd f} \right)  \\ 
- \Re \left( \frac{\derr \fd P \derl f}{\fd f} \right)
 & J_R \end{array}  \right). \ee
The first fundamental form of the surface $\cf$ is
\be\label{1stff}
I  =   J_L {\rm d}\xi_L^2 -  2 \Re \left( \frac{\derr \fd P \derl f}{\fd f} \right)
{\rm d}\xi_L {\rm d}\xi_R + J_R {\rm d}\xi_R^2 .
\ee

In order to establish conditions on a solution $f$ of the Euler--Lagrange equations (\ref{eqnmotN}) under which 
the surface exists, we employ the Schwarz inequality 
\be\label{Schw}
 |\langle a, Ab \rangle|^2 \leq \langle a, Aa \rangle \langle b, Ab \rangle 
\ee
valid for any positive hermitean operator $A$ (see e.g. \cite{BEH}).
Also note that equality in (\ref{Schw}) holds only if there exists $\alpha \in \C$ such that either
$ \langle \alpha a +b, A (\alpha a +b) \rangle =0$ or $\langle a +\alpha  b, A ( a + \alpha b) \rangle=0$ holds.
We may write
\be\label{jdse}
 J_D =\frac{\langle \partial_D f, P \partial_D f \rangle}{\langle f, f \rangle} \geq 0, \; D=L,R 
\ee
and 
\be\label{dgse} 
\det G = \frac{ \langle \derl f, P \derl f \rangle \langle \derr f, P \derr f \rangle - \left( \Re \langle \derl f,
 P \derr f \rangle \right)^2 }{\langle f, f \rangle^2} \geq 0
\ee
since
$$ \langle \derl f, P \derl f \rangle \langle\derr f, P \derr f \rangle \geq | \langle \derl f, P \derr f \rangle|^2 \geq  \left( \Re 
\langle \derl f, P \derr f \rangle \right)^2.$$
Therefore the first fundamental form $I$ defined by (\ref{1stff}) is positive for any solution $f$ 
of the Euler--Lagrange equations (\ref{eqnmotN}).

Analyzing the cases when equalities in (\ref{jdse}),(\ref{dgse}) hold we find that $I$ is 
positive definite in the point $(\xi^0_L,\xi^0_R)$ either if the inequality
\be\label{posdef1} 
\Im \left( \derl \fd(\xi^0_L,\xi^0_R) P \derr f(\xi^0_L,\xi^0_R) \right) \neq 0 
\ee
holds or if the vectors
\be\label{posdef2} 
\derl f(\xi^0_L,\xi^0_R), \derr f(\xi^0_L,\xi^0_R), f(\xi^0_L,\xi^0_R)
\ee
are linearly independent. Therefore any of the conditions (\ref{posdef1}),(\ref{posdef2}) is a sufficient condition for the existence 
of the surface $\cf$ associated with the solution $f$ of the Euler--Lagrange equations 
(\ref{eqnmotN}) in the vicinity of the point $(\xi^0_L,\xi^0_R)$. 
If neither of the conditions (\ref{posdef1}),(\ref{posdef2}) is met on an image ${\rm Im}_X(\Theta)$ 
of a lower--dimensional subset $\Theta\subset\Omega$ then the surface 
$\cf$ may or may not exist, depending on circumstances. If both conditions (\ref{posdef1}),(\ref{posdef2})
are violated in the whole neighborhood $\Upsilon \subset \Omega$ of the point $(\xi^0_L,\xi^0_R)$ 
then the surface doesn't exist in this neighborhood $\Upsilon$.

Using (\ref{metric}) we can write the formula for Gaussian curvature as
\be
K = \frac{1}{\sqrt{J_{L} J_{R} - G_{LR}^2}} \derr 
\left( \frac{\derl G_{LR} - \frac{1}{2} G_{LR} \derl(\ln J_L)}{\sqrt{J_{L} J_{R} - G_{LR}^2}}\right),
\ee
where $$ G_{LR} = - \Re \left( \frac{\derr \fd P \derl f}{\fd f} \right).$$
In the $\C P^1$ case a surprising simplification occcurs and we find that the Gaussian curvature is a negative constant, $ K= -4.$
Consequently, there are no umbilical points on the surface and any regular 
solution of the Euler--Lagrange equations (\ref{eqnmotN}) gives rise to a pseudosphere immersed in $su(2) \simeq \R^3$. Several examples 
of such pseudospheres were present in \cite{Grsnosym}, one is also reproduced in Section \ref{example}.

\section{The Gauss--Weingarten equations}

Now we may formally determine a moving frame on the surface $\cf$ and write 
the Gauss--Weingarten equations in the $\C P^{N-1}$ case. 
Let $f$ be a solution of the Euler--Lagrange equations (\ref{eqnmotN})
 such that ${\rm det}(G)$ is not zero in a neighborhood of a regular point $(\xi_L^0,\xi_R^0)$ in $\Omega$.
Assume also that the surface $\cf$ (\ref{surfaceN}), associated with these equations is described by the moving frame 
$$\vec \tau=(\derl X, \derr X, n_{3}, \ldots, n_{N^2-1})^T,$$
where the vectors $\derl X, \derr X, n_{3}, \ldots, n_{N^2-1}$ satisfy the normalization conditions
$$ (\derl X,\derl X)=J_{L}, \, (\derl X,\derr X)= G_{LR}, \, (\derr X,\derr X)= J_{R}, $$
\be\label{normnorm} (\derl X, n_{k})=(\derr X, n_{k})=0, \, (n_j, n_{k})=\delta_{jk}. \ee
We now show that the moving frame satisfies the Gauss--Weingarten equations
\begin{eqnarray}
\nonumber \derl \derl X & = & A^L_L \derl X + A^L_R \derr X + Q^L_j n_j, \\
\nonumber \derl \derr X & = & \tilde H_j n_j, \\ 
\nonumber \derl n_j & = & \alpha^L_j \derl X + \beta^L_j \derr X +s^L_{jk} n_k, \\
\nonumber \derr \derl X & = & \tilde H_j n_j, \\ 
\nonumber \derr \derr X & = & A^R_L \derl X + A^R_R \derr X + Q^R_j n_j, \\
\derr n_j & = & \alpha^R_j \derl X + \beta^R_j \derr X +s^R_{jk} n_k \label{gweqN},
\end{eqnarray}
where $s^L_{jk}+s^L_{kj}=0$, $s^R_{jk}+s^R_{kj}=0,$ $j,k=3,\ldots, N^2-1$,
$$ \alpha^L_j =  \frac{\tilde H_j G_{LR}- Q^L_j J_{R}}{{\rm det}G}, \, \, \, 
\beta^L_j =  \frac{Q^L_j G_{LR} -  \tilde H_j J_{L}}{{\rm det}G},$$ 
$$ \alpha^R_j =  \frac{Q^R_j G_{LR} -  \tilde H_j J_{R}}{{\rm det}G}, \, \,  \, 
\beta^R_j =  \frac{\tilde H_j G_{LR}-Q^R_j J_{L}}{{\rm det}G},$$ 
\begin{eqnarray}
\nonumber A^L_L  & = & \frac{1}{{\rm det}G} \Re \left\{ 
\frac{1}{\fd f} \left( J_R \derl \fd + G_{LR} \derr \fd \right) P \derl \derl f \right. \\
\nonumber & - & \left. \frac{2\derl \fd f}{(\fd f)^2}  (\derl \fd P \derr f) G_{LR} - \frac{2 \fd \derl f}{\fd f} J_L J_R \right\}, \\
\nonumber A^L_R  & = & \frac{1}{{\rm det}G} \Re \left\{ 
- \frac{1}{\fd f} \left( J_L \derr \fd + G_{LR} \derl \fd \right) P \derl \derl f \right. \\
\label{As} & + & \left. \frac{2\derl \fd f}{(\fd f)^2}  (\derl \fd P \derr f) J_L + \frac{2 \fd \derl f}{\fd f} J_L G_{LR}  \right\},
\end{eqnarray}
and $A^R_L,A^R_R$ have similar form which can be obtained by exchange $L \leftrightarrow R$. 
The explicit form of the coefficients $\tilde H_j,Q^D_j$ (where $D=L,R$; $j=3,\ldots, N^2-1$) depends on 
the chosen orthonormal basis $\{ n_3, \ldots, n_{N^2-1} \}$ of the normal space to the surface $\cf$ 
at the point $X(\xi_L^0,\xi_R^0)$.

Indeed, if $\derl X, \derr X$ are defined by (\ref{dlxdrx}) for an arbitrary solution $f$ of the Euler--Lagrange equations (\ref{eqnmotN}),
then by straightforward calculation using (\ref{eqnmotN}) one finds that
\begin{eqnarray}
\nonumber \derl \derr X  & = &  \derr \derl X = [\derl P, \derr P] = \\
\nonumber & = & \frac{1}{\fd f} \left( P \derl f \otimes \derr \fd P - P \derr f \otimes \derl \fd P \right) \\
 & + & \frac{1}{(\fd f)^2} \left( \derl \fd P \derr f - \derr \fd P \derl f \right) f \otimes \fd. \label{2ndderlr}
\end{eqnarray}
By computing 
\be\label{11}
\tr \left( \derl \derr X . \partial_{D} X \right) = \pm \tr ([\derl P,\derr P].[ \partial_{D} P,P]) = 0, \; D=L,R
\ee
we conclude that $\derl \derr X$ is perpendicular to the surface $\cf$ and consequently it has the form given in (\ref{gweqN}).

The remaining relations in (\ref{gweqN}) and (\ref{As}) follow as differential consequences 
from the assumed normalizations of the normals (\ref{normnorm}), e.g.
$$ (n_j,n_k)=0, \; j \neq k $$
which gives 
$$0 = (\derl n_j,n_k)+(\derl n_k,n_j) = s^L_{jk} + s^L_{kj}.$$
Similarly
$$ (n_j,\derl X)=0, \, \, (n_j,\derr X)=0   $$
by differentiation leads to
$$ (\derr n_j,\derl X) + (n_j,\derl \derr X) =0, \, \, (\derr n_j,\derr X) + (n_j, \derr \derr X)=0 $$
implying
$$ J_{L} \alpha^R_j + G_{LR} \beta^R_j + \tilde H_j =0, \, \, G_{LR} \alpha^R_j + J_{R} \beta^R_j + Q^R_j =0. $$
Consequently, $\alpha^R_j,\beta^R_j$ can be determined in terms of $\tilde H_j, Q^R_j$ and of the components of the induced metric 
$G$. The remaining coefficients $\alpha^L_j,\beta^L_j$ are derived in an analogous way by exchanging indices $\L \leftrightarrow R$
in the successive differentiations.

The coefficients $A^L_L,\ldots,A^R_R$ are obtained by requiring that $(\partial_{D} \partial_{D} X - A^D_L \derl X - A^D_R \derr X)$
is normal to the surface, i.e. 
\be\label{trdx}
\tr \left( \partial_{B} X.(\partial_{D} \partial_{D} X - A^D_L \derl X - A^D_R \derr X) \right)=0, \; B,D =L,R.
\ee
From (\ref{explM}) and (\ref{cldx}) we find 
\begin{eqnarray}
\nonumber \derl \derl X & = & \frac{1}{\fd f} \left( P \derl \derl f \otimes \fd - f \otimes \derl \derl \fd P \right) \\
\nonumber & + & \frac{2}{(\fd f)^2} \left( (\derl \fd f) f \otimes \derl \fd P - (\fd \derl f) P \derl f \otimes \fd   \right), \\
\nonumber \derr \derr X & = & \frac{1}{\fd f} \left( f \otimes \derr \derr \fd P - P \derr \derr f \otimes \fd  \right) \\
 & + & \frac{2}{(\fd f)^2} \left( (\fd \derr f) P \derr f \otimes \fd - (\derr \fd f) f \otimes \derr \fd P \right), \label{2ndder}
\end{eqnarray}
and after substituting the above expressions into (\ref{trdx}) we solve the resulting linear equations for $A^D_B$.

Let us note that the Gauss--Weingarten equations (\ref{gweqN}) can be written equivalently in the $N \times N$ matrix form
\be
\derl \vec \tau = U \vec \tau, \qquad \derr \vec \tau = V \vec \tau,
\ee
where
\begin{eqnarray}
U & = & \left( \begin{array}{ccccc} A^L_L & A^L_R & Q_3^L & \ldots & Q^L_{N^2-1} \\
0 & 0 & \tilde H_3 & \ldots & \tilde H_{N^2-1} \\
\alpha^L_3 & \beta^L_3 & s^L_{33} & \ldots & s^L_{3(N^2-1)} \\
\ldots & \ldots & \ldots & \ldots & \ldots \\
\alpha^L_{(N^2-1)} & \beta^L_{(N^2-1)} & s^L_{(N^2-1)3} & \ldots & s^L_{(N^2-1)(N^2-1)} 
\end{array}
\right), \nonumber \\
V & = & \left( \begin{array}{ccccc} 0 & 0 & \tilde H_3 & \ldots & \tilde H_{N^2-1} \\
A^R_L & A^R_R & Q_3^R & \ldots & Q^R_{N^2-1} \\
\alpha^R_3 & \beta^R_3 & s^R_{33} & \ldots & s^R_{3(N^2-1)} \\
\ldots & \ldots & \ldots & \ldots & \ldots \\
\alpha^R_{(N^2-1)} & \beta^R_{(N^2-1)} & s^R_{(N^2-1)3} & \ldots & s^R_{(N^2-1)(N^2-1)} 
\end{array} \right).
\end{eqnarray}

The Gauss--Codazzi--Ricci equations 
\be\label{gcod1}
\derr U - \derl V + [U,V] =0
\ee
are compatibility conditions for the Gauss--Weingarten equations (\ref{gweqN}).
They are the necessary and sufficient conditions for the local existence of the corresponding surface $\cf$.
It can be easily checked that they are identically satisfied for any solution $f$ of the Euler--Lagrange
equations (\ref{eqnmotN}).

The second fundamental form of the surface $\cf$ at the regular point $p$ 
takes in general\footnote{In the familiar $\R^3$ case 
the normal space $N_p \cf$ is conventionally identified with $\R$.} the shape of a map
$$ {\bf II}(p): T_p \cf \times T_p \cf \rightarrow N_p \cf, $$
where $T_p \cf$, $N_p \cf$ denote the tangent and normal space to the surface $\cf$ at the point $p$, respectively.
According to \cite{Kob,Wil}, the second fundamental form and the mean curvature vector can be expressed as
\bea
\label{IIN1} {\bf II}  =  (\derl \derl X)^{\perp} {\rm d} \xi_L {\rm d} \xi_L + 2 (\derl \derr X)^{\perp} {\rm d} \xi_L {\rm d} \xi_R
+(\derr \derr X)^{\perp} {\rm d} \xi_R {\rm d} \xi_R, & & \\
\label{HN1}
{\bf H}  =  \frac{1}{\det G} \left( J_{R} (\derl \derl X)^{\perp}  - 2 G_{LR} (\derl \derr X)^{\perp} 
+ J_{L} (\derr \derr X)^{\perp} \right), & &
\eea
where $(\;)^{\perp}$ denotes the normal part of the vector.  
In our case, given the decomposition of $\partial_D \partial_B X$ into the tangent and normal parts in the Gauss--Weingarten
equations (\ref{gweqN}), the expressions (\ref{IIN1}),(\ref{HN1}) take the form
\begin{eqnarray}
\nonumber {\bf II} & = & (\derl \derl X-A^L_L \derl X - A^L_R \derr X) {\rm d} \xi_L {\rm d} \xi_L 
+ 2 (\derl \derr X) {\rm d} \xi_L {\rm d} \xi_R +\\
& + &  (\derr \derr X-A^R_L \derl X - A^R_R \derr X) {\rm d} \xi_R {\rm d} \xi_R, \label{IIN2} \\
\nonumber {\bf H} & = & \frac{1}{\det G} \left( \right. J_{R} (\derl \derl X-A^L_L \derl X - A^L_R \derr X) 
- 2 G_{LR}  (\derl \derr X)  + \\
& + &  J_{L} (\derr \derr X-A^R_L \derl X - A^R_R \derr X)  \left. \right). \label{HN2}
\end{eqnarray}
Consequently, the Willmore functional \cite{Wil} is
\be
W = \int |{\bf H}|^2 \sqrt{{\rm det G} } {\rm d} \xi_L {\rm d} \xi_R.
\ee
The derivatives $\partial_D \partial_B X$ are expressed explicitly in terms of $f$ in equations (\ref{2ndderlr}) and(\ref{2ndder}).
Unfortunately, it is clear that after explicit calculation of $(\partial_B \partial_D X)^{\perp}$ 
in the case of $N>2$, both the second fundamental form and the mean curvature vector contain terms like 
$ P \derl \derl f \otimes \fd $ etc., which are neither cancelled out by other terms nor projected out by the normal projection.
Therefore the resulting expressions are rather complicated and, for lack of space, we do not present them here.

In the $\C P^{1}$ case the formulae (\ref{IIN2}) and (\ref{HN2}) simplify to
$$ {\bf II} = - 2 (\derr \fd P \derl f  - \derl \fd P \derr f )\ (\1-2 P) \ {\rm d}\xi_L {\rm d}\xi_R, $$
$$ {\bf H}  =  2 \frac{\derr \fd P \derl f  + \derl \fd P \derr f }{\derr \fd P \derl f  - \derr \fd P \derr f } \ (\1-2 P),$$
where the normal to the surface $\cf$ is given by 
\be
n= i (\1 - 2P) \in su(2).
\ee

\section{The moving frame of a surface in the algebra $su(N)$}

Now we proceed to construct the moving frame of the surface $\cf$ immersed in $su(N)$ algebra, i.e. matrices 
$ \derl X,\derr X, n_a, \ a=3,\ldots,N^2-1$ satisfying (\ref{normnorm}). 

Let $f$ be a solution of the Euler--Lagrange equations (\ref{eqnmotN}) and let $(\xi_L^0,\xi_R^0)$ be a regular point 
in $\Omega$, i.e. such that ${\rm det}G(f(\xi_L^0,\xi_R^0)) \neq 0$. Let us denote $f^0=f(\xi_L^0,\xi_R^0)$, 
$X^0=X(\xi_L^0,\xi_R^0)$. Taking into account that
$$ \tr(A)= \tr(\Phi A \Phi^{\dagger}), \; \;  A \in su(N), \  \Phi \in SU(N), $$
we employ the adjoint representation of the group $SU(N)$ in order to bring $ \derl X,\derr X, n_a$ to the simplest form possible. 
We shall request
\be\label{Phireq} \Phi^\dagger f^0  =  (\sqrt{\fdo f^0},0,\ldots,0)^{T}. \ee
Let us choose an orthonormal basis in $su(N)$ in the following form
\begin{eqnarray}
\nonumber  (A_{jk})_{ab} & = &  i (\delta_{ja} \delta_{kb} + \delta_{jb} \delta_{ka} ), \; \; 1\leq j<k\leq N, \\
\nonumber  (B_{jk})_{ab} & = &   (\delta_{ja} \delta_{kb} - \delta_{jb} \delta_{ka} ), \; \; 1\leq j<k\leq N, \\
(C_{p})_{ab} & = & i \sqrt{\frac{2}{p(p+1)}} \left( \sum_{d=1}^p \delta_{da} \delta_{db} - p \delta_{p+1,a} \delta_{p+1,b} 
\right), \; \;  1 \leq p \leq N-1.
\end{eqnarray}

The existence of $\Phi \in SU(N)$ satisfying (\ref{Phireq}) follows from the fact that the $SU(N)$ group acts transitively 
 on the set $\{ a \in \C^N, a^\dagger a = \alpha \}$, where $\alpha \in \R^{+}$.
It should be noted that such $\Phi$ is not unique. A concrete form of $\Phi$ can be constructed as follows:
starting from a general element $a = (a_1,\ldots,a_N)^T  $ 
of $\C^N$ one firstly finds a transformation $\Phi^\dagger_{N-1}$ which transforms $a$ into  
the vector 
$$  a^{(N-1)} = (a_1,\ldots,a_{N-2},\sqrt{ a_{N-1} \bar{a}_{N-1} + a_{N} \bar{a}_{N} },0), \; \; a^{(N-1)} \bar{a}^{(N-1)}= |a|^2.$$
It is easy to see that the desired transformation is
$$ 
 \Phi^\dagger_{N-1}   =  \left( \begin{array}{cccc} 
1 & 0 &  \ldots &  0 \\
0 & 1 &  \ldots &  0 \\
\vdots & \vdots & \vdots & \vdots  \\
0 & \ldots  & \frac{{\bar{a}_{N-1}}}{ \left(a_{N-1} \bar{a}_{N-1} + a_{N} \bar{a}_{N} \right)^{1/2} }  & 
\frac{{\bar{a}_{N}}}{\left(a_{N-1} \bar{a}_{N-1} + a_{N} \bar{a}_{N}\right)^{1/2}}  \\
0 & \ldots  & 
- \frac{{a_{N} }}{\left(a_{N-1} \bar{a}_{N-1} + a_{N} \bar{a}_{N}\right)^{1/2}} 
& \frac{{a_{N-1}}}{\left(a_{N-1} \bar{a}_{N-1} + a_{N} \bar{a}_{N}\right)^{1/2}}  
\end{array}  \right).$$
In the next step one eliminates the last but one element of $a^{(N-1)}$ using the transformation
$$  \Phi^\dagger_{N-2}   =  \left( \begin{array}{ccccc} 
1 & 0 &  0 &  \ldots &  0 \\
0 & 1 &  0 &  \ldots &  0 \\
\vdots & \vdots & \vdots & \vdots & \vdots  \\
0 & \ldots  & \frac{{\bar{a}_{N-2}}}{\left( \sum_{j=N-2}^{N} a_j \bar{a}_j \right)^{1/2}}  & 
\frac{\left( \sum_{j=N-1}^{N} a_j  \bar{a}_j\right)^{1/2}}{\left( \sum_{j=N-2}^{N} a_j \bar{a}_j \right)^{1/2}}   & 0   \\
0 & \ldots  & 
- \frac{\left( \sum_{j=N-2}^{N} a_j \bar{a}_j \right)^{1/2}}{\left( \sum_{j=k}^{N} a_j \bar{a}_j \right)^{1/2}} 
& \frac{{a_{N-2}}}{\left( \sum_{j=N-2}^{N} a_j \bar{a}_j \right)^{1/2}}  & 0  \\
0   & 0 & \ldots & 0 & 1 \end{array}  \right).$$
This gives	
$$ a^{(N-2)} = \Phi^\dagger_{N-2} a^{(N-1)} = \Phi^\dagger_{N-2} \Phi^\dagger_{N-1} a = 
(a_1,\ldots,a_{N-3},\left( \sum_{j=N-2}^{N} a_j \bar{a}_j  \right)^{1/2},0,0). $$
By induction and redesignation $a=f^0$ one arrives at the explicit form of $\Phi$ 
\begin{eqnarray}
\nonumber \Phi^\dagger_k  & = & \left( \begin{array}{cccccc} 
1 & 0 &  \ldots &  0 & \ldots &  0 \\
0 & 1 &  \ldots &  0 & \ldots &  0 \\
\vdots &  &  & \vdots &  & \vdots \\
0 & \ldots  & \frac{{\bar{f}_{k}^0}}{\left( \sum_{j=k}^{N} f_j^0 \bar{f}_j^0 \right)^{1/2}}  & 
\frac{\left( \sum_{j=k+1}^{N} f_j^0  \bar{f}_j^0\right)^{1/2}}{\left( \sum_{j=k}^{N} f_j^0 \bar{f}_j^0 \right)^{1/2}}   & 0 & \ldots  \\
0 & \ldots  & 
- \frac{\left( \sum_{j=k+1}^{N} f_j^0 \bar{f}_j^0 \right)^{1/2}}{\left( \sum_{j=k}^{N} f_j^0 \bar{f}_j^0 \right)^{1/2}} 
& \frac{{f_{k}^0}}{\left( \sum_{j=k}^{N} f_j^0 \bar{f}_j^0 \right)^{1/2}}  & 0 & \ldots \\
0  & \ldots & 0 & 0 & 1 & \ldots \\
\vdots &  & \vdots & \vdots &  & \vdots \\
0  & \ldots & 0 & \ldots & 0 & 1
\end{array}  \right)  , \  k \leq N-2, \\
\end{eqnarray}
\begin{eqnarray}
\nonumber \Phi^\dagger_{N-1}  & = & \left( \begin{array}{cccc} 
1 & 0 &  \ldots &  0 \\
0 & 1 &  \ldots &  0 \\
\vdots & \vdots & \vdots & \vdots  \\
0 & \ldots  & \frac{{\bar{f}_{N-1}^0}}{ \left(f^0_{N-1} \bar{f}^{0}_{N-1} + f^0_{N} \bar{f}^{0}_{N} \right)^{1/2} }  & 
\frac{{\bar{f}_{N}^0}}{\left(f^0_{N-1} \bar{f}^{0}_{N-1} + f^0_{N} \bar{f}^{0}_{N}\right)^{1/2}}  \\
0 & \ldots  & 
- \frac{{f_{N}^0 }}{\left(f^0_{N-1} \bar{f}^{0}_{N-1} + f^0_{N} \bar{f}^{0}_{N}\right)^{1/2}} 
& \frac{{f_{N-1}^0}}{\left(f^0_{N-1} \bar{f}^{0}_{N-1} + f^0_{N} \bar{f}^{0}_{N}\right)^{1/2}}  
\end{array}  \right), 
\\ \Phi & = & \Phi_{N-1} \Phi_{N-2}\cdots \Phi_{1}, \; \Phi^\dagger  =  \Phi_{1}^{\dagger} \Phi_{2}^{\dagger} \cdots \Phi_{N-1}^{\dagger} \ \in SU(N) \label{PhiN}.
\end{eqnarray}
If any of the denominators vanishes then the corresponding matrix $\Phi_{k}$ is defined to be the unit matrix.
It is also clear that the group element $\Phi$ constructed in this way is a smooth function of $f,\fd$ and consequently of $\xi_L,\xi_R$.
Thus we find
\begin{eqnarray}
\nonumber \Phi^\dagger f^0 & = & (\sqrt{\fdo f^0},0,\ldots,0)^{T}, \\
\nonumber \derl^\Phi X^0 & \equiv & \Phi^\dagger \derl X(\xi_L^0,\xi_R^0) \Phi  =    \frac{1}{\sqrt{\fdo f^0}}  \left( \begin{array}{cc} 0 & - {\derl^{\Phi} f^0}^\dagger \\ 
{\derl^{\Phi} f^0} & {\bf 0} \end{array} \right), \\
\derr^\Phi X^0 & \equiv & \Phi^\dagger \derr X(\xi_L^0,\xi_R^0) \Phi  =  - \frac{1}{\sqrt{\fdo f^0}} \left( \begin{array}{cc} 0 & - {\derr^{\Phi} f^0}^\dagger \\ 
{\derr^{\Phi} f^0} & {\bf 0} \end{array} \right), \label{movframe0}
\end{eqnarray}
where ${\bf 0}$ denotes the null $(N-1)\times (N-1)$ matrix and the vectors  ${\partial_{D}^{\Phi} f^0}\in \C^{N-1}$ are defined by
$$ ({\partial_{D}^{\Phi} f^0})_{j-1}= (\Phi^\dagger \partial_{D} f(\xi_L^0,\xi_R^0))_{j}, \, D=L,R, \, j=2,\ldots,N.$$

The construction of the moving frame is now straightforward.
Assume that one finds, using a variant of Gramm-Schmidt orthogonalization procedure, the orthonormal vectors 
$$ \tilde A_{1j}, \tilde B_{1j}, \; j=3,\ldots, N$$
satisfying
$$ (\partial_{D}^\Phi X^0, \tilde A_{1j}) =0, \, (\partial_{D}^\Phi X^0, \tilde B_{1j}) =0 $$
and 
\be\label{GrammSch}
 {\rm span} (\partial_{D}^\Phi X^0,\tilde A_{1j}, \tilde B_{1j})_{D=L,R, \ j=3,\ldots,N} = 
{\rm span} (A_{1j},B_{1j})_{j=2,\ldots,N}.
\ee
We identify the remaining tilded and untilded matrices
$$ \tilde A_{jk} = A_{jk}, \ \tilde B_{jk} = B_{jk}, \ \tilde C_{p} = C_{p}, \ 1 < j < k \leq N, \ 1 \leq p \leq N-1. $$
From (\ref{movframe0}) directly follows that
$$ (\partial_{D}^\Phi X^0, \tilde A_{jk})=(\partial_{D}^\Phi X^0,\tilde B_{jk})=(\partial_{D}^\Phi X^0,\tilde C_{p})=0, 
\;  \ 1<j<k\leq N,p<N $$
and as a result of Gramm-Schmidt orthogonalization we get
$$ (\partial_{D}^\Phi X^0,\tilde A_{1k})=(\partial_{D}^\Phi X^0,\tilde B_{1k})=0 $$
and for $1 < j < k \leq N, \ p=1,\ldots, N-1, \  \ i=3,\ldots,N$
$$ (\tilde A_{1i},\tilde A_{jk}) = (\tilde A_{1i},\tilde B_{jk}) = (\tilde A_{1i},\tilde C_{p}) = 
(\tilde B_{1i},\tilde A_{jk}) = (\tilde B_{1i},\tilde B_{jk}) = (\tilde B_{1i},\tilde C_{p}) = 0 . $$
Therefore, under the above given assumptions and notation, we can state the following 
\begin{prop}
The moving frame of the surface $\cf$ at the point $X^0=X(\xi_L^0,\xi_R^0)$ 
\begin{eqnarray}
\nonumber \derl X  & = & \Phi \derl^\Phi X^0 \Phi^\dagger, \\ 
\nonumber \derr X  & = & \Phi \derr^\Phi X^0 \Phi^\dagger, \\
\nonumber n^A_{jk} & = & \Phi \tilde A_{jk} \Phi^\dagger, \\
\nonumber n^B_{jk} & = & \Phi \tilde B_{jk} \Phi^\dagger, \ 1 \leq j < k \leq N, \\
 n^C_{p} & = & \Phi \tilde C_{p} \Phi^\dagger, \ 1 \leq p \leq N-1. \label{movframeN}
\end{eqnarray}
satisfies the normalization conditions (\ref{normnorm})
and consequently the Gauss--Weingarten equations (\ref{gweqN}).
\end{prop}
\smallskip

Note that the first two lines of (\ref{movframeN}) are equivalent to (\ref{movframe0}). 
The remaining lines of (\ref{movframeN}) give a rather explicit 
description of normals to the surface $\cf$. In the $\C P^1$ case a significant simplification occurs, namely there is 
only one normal vector $n^C_1= i \Phi \sigma_3 \Phi^{-1} $ to the surface immersed in $su(2)$ and no orthogonalization is needed.

In the case of $N>2$ the explicit form of the moving frame (\ref{movframeN}) might be quite complicated
 because of the orthogonalization process involved in the construction of 
$$  n^A_{1j},  n^B_{1j}, \; j=3,\ldots, N$$
(i.e. in the construction of $\tilde A_{1j}, \tilde B_{1j}, \; j=3,\ldots, N$).
On the other hand, the remaining normals 
$$  n^A_{jk}, \  n^B_{jk}, \  n^C_{p}, \ 1 < j < k \leq N, \ 1 \leq p \leq N-1, $$
can be constructed without any difficulty. If we chose other group element $\Phi$ satisfying (\ref{Phireq}),
the constructed normals would have been rotated by a local (gauge) transformation from the subgroup of $SU(N)$ leaving 
$\derl X(\xi_L^0,\xi_R^0),\derr X(\xi_L^0,\xi_R^0)$ invariant.

It is also worth noting that from the equations (\ref{2ndderlr}),(\ref{2ndder}) immediately follows that 
\begin{eqnarray}
\nonumber \derl \derl X, \derr \derr X \in \ \Phi \ {\rm span} (A_{1j},B_{1j})_{j=2,\ldots,N} \ \Phi^{\dagger}, \\
\derl \derr X \in \ \Phi \ {\rm span} (A_{jk},B_{jk},C_{p})_{1<j<k\leq N,p<N} \ \Phi^{\dagger},
\end{eqnarray}
i.e.
\begin{eqnarray}
\nonumber (\derl \derl X)^\perp, (\derr \derr X)^\perp \in \  {\rm span} (n^A_{1j},n^B_{1j})_{j=3,\ldots,N}, \\
(\derl \derr X)^\perp=\derl \derr X \in \  {\rm span} (n^A_{jk},n^B_{jk},n^C_{p})_{1<j<k\leq N,p<N}.
\end{eqnarray}

Concerning other possible constructions of the normals, one can observe that one may construct immediately two unit 
normals\footnote{Note that the scalar product
$$(n_P,n_{[\derl X, \derr X]}) \neq 0$$
so that their orthogonalization would be needed.}
$$ n_P = i \sqrt{2} \left( \sqrt{\frac{N-1}{N}} \1 - \sqrt{\frac{N}{N-1}} P \right), $$
$$ n_{[\derl X, \derr X]} = \frac{i[\derl X, \derr X ]}{|[\derl X, \derr X ] |}. $$
In the $su(2)$ case the normals $n_P,n_{[\derl X, \derr X]},n^C_1$ coincide up to the choice of orientation,
but in general the relation of $ n_P, n_{[\derl X, \derr X]}$ 
to $n^A_{jk},n^B_{jk},n^C_{p}$ is rather complicated and difficult to express in a closed form. 
In principle one could attempt to construct the moving frame directly from these normals
by taking normal parts of commutators of them with $\derl X, \derr X$ etc.\footnote{This 
can be proved to be possible at least in the $su(3)$ case by observing that $\derl X, \derr X$ 
generate via commutators the whole algebra $su(3)$.},
without need to construct the group element  $\Phi$. Unfortunately, such procedure does not seem to be 
computationally feasible at the moment, leaving this subject open for further investigation.

\section{Example of surface in the algebra $su(2)$}\label{example}

\begin{figure}[thb]
\epsfxsize=5in
\begin{center}       
\leavevmode 
\epsffile{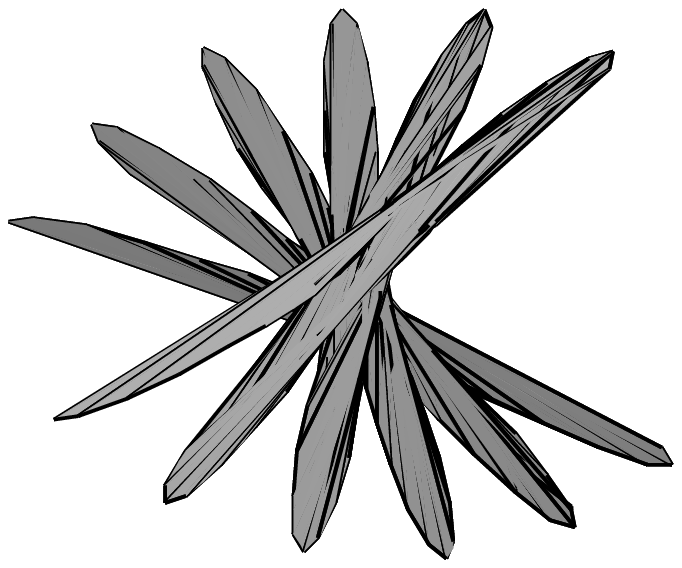}
\end{center}
\caption{The surface associated with the solution (\ref{expwell}), $p=-3/2$} 
\label{ellsolfig}
\end{figure}

As an example of a surface obtained using the described method we present a picture of a surface in $su(2)$
associated with the $\C P^1$ sigma model on Minkowski space (see Figure \ref{ellsolfig}). 
The following solution of the Euler--Lagrange equations (\ref{eqnmotN})
in this case was obtained by us using the symmetry reduction method
\be\label{expwell}
f = \left( 1, \sqrt{\frac{(p-1) \cosh(g(\chi))+(p+1)}{(p-1) \cosh(g(\chi))-(p+1)} } \ 
{\rm exp} ( i(\xi_L-h(\chi)) )\right),
\ee
where
\begin{eqnarray}
\nonumber h(\chi) & =  &
\arctan \left( \frac{p+1}{2\sqrt{-p}} \tanh g(\chi) \right) +\frac{( p + 2 \sqrt{-p} -1) \chi }{2(p-1)}, \\
\nonumber g(\chi) & = &\frac{(p+1) \chi }{2(p-1)},   \; \chi= \xi_L-\xi_R, \; p<-1.
\end{eqnarray}
The formulae for the first and second fundamental forms, moving frame etc. of the associated  surface $\cf$ 
were obtained but are too lengthy to be presented here. The computation of the surface, i.e. the Weierstrass 
representation (\ref{surfaceN}),
was performed numerically. The Gaussian curvature is $K=-4$, the mean curvature is
$$H = - \frac{{\rm e}^{4 g(\chi)}-6 {\rm e}^{2 g(\chi)} +1}{2 \ {\rm e}^{g(\chi)} \ ({\rm e}^{2 g(\chi)}-1)}.$$

Other examples were presented in \cite{Grsnosym}.

\section{Conclusions}\label{concl}

The main purpose of this paper was to provide the structural equations of 2--dimensional orientable 
smooth surfaces immersed in $su(N)$ algebra. 
The surfaces were obtained from the $\C P^{N-1}$ sigma model defined on 2--dimensional Minkowski space.

The most important advantage of the method presented is that it gives effective tools for constructing surfaces without 
reference to additional considerations, proceeding directly from the given $\C P^{N-1}$ model equations (\ref{eqnmotN}). 
We demonstrated through the use of Cartan's language of moving frames that one can derive via $\C P^{N-1}$ model the first 
and second fundamental forms of the surface as well as the relations between them as expressed in the Gauss--Weingarten and 
Gauss--Codazzi--Ricci equations. 
We presented an extension of the classical Enneper--Weierstrass representation of surfaces in multi--dimensional spaces,
expressed in terms of any nonsingular (i.e. such that $\det G \neq 0$)  solution of the $\C P^{N-1}$ sigma model.
In particular, we showed that in the $\C P^1$ case such description of surfaces in the $su(2)$ algebra leads to 
constant negative Gaussian curvature surfaces.

Our method can be particularly useful in applications like the theory of phase transitions or fluid membranes \cite{Nel,Cha,Ou,Saf}, 
where numerical
aproaches have prevailed so far. Even in cases when the Weierstrass representation of a surface cannot be integrated explicitly,
 the surface's main characteristics can be derived in analytical form which lends itself to physical interpretations. 

\subsection*{Acknowledgments}

This work was supported in part by research grants from NSERC of Canada.
Libor \v Snobl acknowledges a postdoctoral fellowship awarded by the Laboratory of Mathematical Physics of the CRM, 
Universit\'e de Montr\'eal. The authors thank Pavel Winternitz for helpful and interesting discussions on the topic
of this paper.

\end{document}